\documentclass{article} 

\usepackage{amsmath, amssymb}

\title{THE TRACE FORMULA WITH RESPECT TO THE TWISTED GROVER MATRIX OF A MIXED DIGRAPH}
\author{Takashi KOMATSU \\
Department of Bioengineering, School of Engineering,\\ 
The University of Tokyo \\
Bunkyo, Tokyo, 113-8656, JAPAN \\ 
e-mail: komatsu@coi.t.u-tokyo.ac.jp \\ 
Sho KUBOTA\\
Department of Applied Mathematics, 
Faculty of Engineering, \\ 
Yokohama National University, \\
Hodogaya, Yokohama 240-8501, JAPAN \\
e-mail: kubota-sho-bp@ynu.ac.jp \\
Norio KONNO\\
Department of Applied Mathematics, 
Faculty of Engineering, \\ 
Yokohama National University, \\
Hodogaya, Yokohama 240-8501, JAPAN \\
e-mail: konno-norio-bt@ynu.ac.jp \\
Iwao SATO \\ 
National Institute of Technology, Oyama College, \\ 
Oyama, Tochigi 323-0806, JAPAN \\
e-mail: isato@oyama-ct.ac.jp }

\begin{document}
 \maketitle

\vspace{5mm}

Running head tiltle:

Trace formula for the twisted Grover matrix 

\vspace{5mm}

The address for manuscript correspondence: 

Iwao Sato 

Oyama National College of Technology, 
Oyama, Tochigi 323-0806, JAPAN

Tel: 0285-20-2176

Fax: 0285-20-2880

E-mail: isato@oyama-ct.ac.jp

\clearpage

\begin{abstract}
We define a zeta function woth respect to the twisted Grover matrix of a mixed digraph, and  
present an exponential expression and a determinant expression of this zeta function. 
As an application, we give a trace formula with respect to the twisted Grover matrix of a mixed digraph. 
\end{abstract}

\vspace{2mm} 


\vspace{2mm}

\clearpage

\section{Introduction}

The study of zeta functions of graphs was started by Ihara [6]. 
In [6], he defined the Ihara zeta functions of graphs, and showed that the 
reciprocals of the Ihara zeta functions of regular graphs are explicit polynomials. 
Hashimoto [5] treated multivariable zeta functions of bipartite graphs. 
Bass [2] generalized Ihara's result on the Ihara zeta function of 
a regular graph to an irregular graph $G$. 

The trace formula for a connected graph $G$ is closely related 
to the Ihara zeta function of $G$. 
Ahumada [1] gave a trace formula for a regular graph(c.f., [17,19]). 
McKay [13] determined the limiting probability density for the eigenvalues of 
a series of regular graphs. 
Sunada [16] presented the semicircle law for the distribution of 
eigenvalues of regular graphs when their girths and degrees are divergent. 

The transition matrix of a discrete-time quantum walk on a graph 
is closely related to the Ihara zeta function of a graph. 
As a quantum counterpart of the classical random walk, a quantum walk has recently 
attracted much attention for various fields, and quantum walks on graphs have been studied 
by many researchers(see Konno [8], Venegas-Andraca [18], Portugal [14]). 
Konno [7] presented the limit theorem (the Konno distribution) of a two-state quantum walk on $\mathbb{Z} $. 
The Konno distribution is quite different from the normal distribution, but is similar to the arcsin law.  
Furthermore, various limit theorems for the probability of a quantum walk on a graph are obtained. 

A discrete-time quantum walk is a quantum process on a graph whose state vector is governed by 
a matrix called the transition matrix.
A relationship between the Grover walk and the Ihara zeta function 
of a graph was given by Ren et al. [15]. 
The Grover walk [4] is a discrete-time quantum walk on a graph which originates from the Grover algorithm. 
The Grover matrix that is the transition matrix of the Grover walk is a typical transition matrix 
of a discrete-time quantum walk on a graph.  
Konno and Sato [9] presented a formula for the characteristic polynomial of the Grover matrix of a graph. 
Konno, Mitsuhashi, Morita and Sato [10] presented a trace formula with respect to the Grover matrix of a regular graph. 
Kubota, Segawa and Taniguchi [11] introduced a twisted Grover matrix of a mixed digraph as a generalization of 
the Grover matrix of a graph, and gave its spectral analysis. 

In this paper, we give a trace formula with respect to the twisted Grover matrix of a regular mixed digraph. 
In Section 2, we give a short review for the Ihara zeta function of a graph, and a trace formula with respect to 
the Ihara zeta function of a regular graph.  
In Section 3, we present a review for a zeta function and a trace formula with respect to the Grover matrix of a graph. 
In Section 4, we introduce a zeta function with respect to the twisted Grover matrix of a mixed digraph, and 
we present its determinant expression. 
Furthermore, we treat its poles in the case of a regular mixed digraph, and present an exponential expression 
and the Euler product of this zeta function. 
In Section 5, we present a trace formula with respect to the twisted Grover matrix of a regular mixed digraph. 

For the Selberg trace formula, the reader is referred to [17,19], respectively. 

\section{Preliminaries}

Graphs and digraphs treated here are finite, simple and unweighted.
Let $G$ be a connected graph and $D$ the symmetric digraph 
corresponding to $G$. 
Set $D(G)= \{ (u,v),(v,u) \mid uv \in E(G) \} $. 
We also refer $D$ as a graph $G$. 
For $e=(u,v) \in D(G)$, $u=o(e)$ and $v=t(e)$ are the {\em origin} and the {\em terminus} of $e$, respectively.  
Furthermore, let $e^{-1}=(v,u)$ be the {\em inverse} of $e=(u,v)$. 

A {\em path $P$ of length $n$} in $D$(or $G$) is a sequence 
$P=(e_1, \ldots ,e_n )$ of $n$ arcs such that $e_i \in D(G)$,
$t( e_i )=o( e_{i+1} )(1 \leq i \leq n-1)$. 
Set $ \mid P \mid =n$, $o(P)=o( e_1 )$ and $t(P)=t( e_n )$. 
Also, $P$ is called an {\em $(o(P),t(P))$-path}. 
We say that a path $P=(e_1, \cdots ,e_n )$ has a {\em backtracking} 
if $ e^{-1}_{i+1} =e_i $ for some $i(1 \leq i \leq n-1)$. 
A $(u,v)$-path is called an {\em $u$-cycle} 
(or {\em $u$-closed path}) if $u=v$. 
The {\em inverse cycle} of a cycle 
$C=( e_1, \ldots ,e_n )$ is the cycle 
$C^{-1} =( e^{-1}_n , \cdots ,e^{-1}_1 )$.

We introduce an equivalence relation between cycles. 
Two cycles $C_1 =(e_1, \ldots ,e_m )$ and 
$C_2 =(f_1, \cdots ,f_m )$ are called {\em equivalent} if 
$f_j =e_{j+k} $ for all $j$, where the sum $j+k$ is performed modulo $m$. 
The inverse cycle of $C$ is in general not equivalent to $C$. 
Let $[C]$ be the equivalence class that contains a cycle $C$. 
A cycle $C=( e_1, \cdots ,e_n )$ has a {\em tail} if $e^{-1}_n =e_1 $.  
A cycle $C$ is {\em reduced} if $C$ has neither a backtracking nor a tail.  
Let $B^r$ be the cycle obtained by going $r$ times around a cycle $B$. 
Such a cycle is called a {\em multiple} of $B$. 
Furthermore, a cycle $C$ is {\em prime} if it is not a multiple of 
a strictly smaller cycle. 
Note that each equivalence class of prime, reduced cycles of a graph $G$ 
corresponds to a unique conjugacy class of 
the fundamental group $ \pi {}_1 (G,v)$ of $G$ at a vertex $v$ of $G$. 
Then the {\em Ihara zeta function} ${\bf Z} (G, u)$ of a graph $G$ is 
defined to be the function of $u \in {\cal C}$ with $|u|$ sufficiently small, 
given by 
\[
{\bf Z} (G,u)= {\bf Z}_G (u)= \prod_{[C]} \  (1- u^{ \mid C \mid } )^{-1}, 
\] 
where $[C]$ runs over all equivalence classes of prime, reduced cycles 
of $G$.

Let $G$ be a connected graph with $n$ vertices $v_1, \ldots ,v_n $, 
and $n \in \mathbb{N} $. 
The {\em adjacency matrix} ${\bf A}= {\bf A} (G)=(a_{ij} )$ is 
the square matrix such that $a_{ij} =1$ if $v_i$ and $v_j$ are adjacent, 
and $a_{ij} =0$ otherwise.
Let $Spec(G)$ be the set of all eigenvalues of ${\bf A} (G)$. 
Let ${\bf D} =( d_{ij} )$ be the diagonal matrix with 
$d_{ii} = \deg {}_G \  v_i $, and ${\bf Q} = {\bf D} -{\bf I} $. 
The {\em degree} $ \deg {}_G \  v= \deg v$ of a vertex $v$ in $G$ is 
defined by $ \deg {}_G \  v = \mid \{ w \mid vw \in E(G) \} \mid $. 
A graph $H$ is called {\em $k$-regular} if $ \deg {}_H v=k$ for 
each vertex $v \in V(H)$.

\newtheorem{theorem}{Theorem}
\begin{theorem}[Ihara]
Let $G$ be a connected $(q+1)$-regular graph with $n$ vertices. 
Set $Spec(G)= \{ \lambda {}_1 , \ldots , \lambda {}_n \} $. 
Then the reciprocal of the Ihara zeta function of $G$ is 
\[
{\bf Z} (G,u)^{-1} = (1- u^2 ) {}^{m-n} 
\det ( {\bf I}_n -u {\bf A} (G)+q u^2 {\bf I}_n ) 
= (1- u^2 ) {}^{(q-1)n/2} 
\prod^{n}_{j=1} (1- \lambda {}_j u+q u^2 )
\]
where $m=|E(G)|$. 
\end{theorem}

Let $G$ be a connected $(q+1)$-regular graph. 
Furthermore, let $h( \theta )$ be a complex-valued function on $\mathbb{R} $ that 
satisfies the following properties: 
\begin{enumerate}
\item $h( \theta + 2 \pi )=h( \theta )$, 
\item $h(- \theta )=h( \theta )$, 
\item $h( \theta )$ is analytically continuable to an analytic function over 
$Im \  \theta < \frac{1}{2} \log q + \epsilon ( \epsilon >0)$. 
\end{enumerate}
For this $h( \theta )$, we define its {\em Fourier transform} by 
\[
\hat{h} (k)= \frac{1}{ 2 \pi } \int^{ 2 \pi }_0 h( \theta ) 
e {}^{ \sqrt{-1} k \theta } d \theta , 
\]
where $k \in \mathbb{Z} $.

\begin{theorem}[Ahumada]
Let $G$ be a connected $(q+1)$-regular graph with $n$ vertices. 
Set $Spec(G)= \{ \lambda {}_1 , \ldots , \lambda {}_n \} $. 
Let $ \lambda {}_1 , \ldots , \lambda {}_l $ be the eigenvalues of $G$ 
for which $1- \lambda {}_j u+ q u^2 =0$ has imaginary roots. 
Furthermore, for each $ \lambda {}_i (1 \leq i \leq l)$, 
let $ q^{-1/2} e {}^{ \sqrt{-1} \theta {}_i}$ be a root of 
$1- \lambda {}_j u+ q u^2 =0$. 
Then the following trace formula holds:
\[
\sum^{l}_{i=1} h( \theta {}_i )= \frac{2nq(q+1)}{ \pi } 
\int^{ \pi }_0 \frac{ \sin {}^2 \theta }{(q+1 )^2 -4q \cos {}^2 \theta } 
h( \theta ) d \theta 
+ \sum_{[C]} \sum^{ \infty }_{m=1} \mid C \mid q^{-m \mid C \mid /2} 
\hat{h} ( m \mid C \mid ) , 
\]
where $[C]$ runs over all equivalence classes of prime, reduced cycles 
of $G$, 
\end{theorem}

Let $G$ be a connected $(q+1)$-regular graph with $n$ vertices. 
Furthermore, let $Spec (G) = \{ \lambda {}_1 , \ldots , \lambda {}_n \} $. 
By Theorem 1, the poles of ${\bf Z} (u)$ are $\pm 1$ and roots of 
$1- \lambda {}_j u+q u^2 =0 \ (1 \leq j \leq n)$. 
Therefore, $u=  q^{-1/2} e {}^{ \sqrt{-1} \theta } $ is a pole of 
${\bf Z} (u)$ if and only if 
$ \lambda = 2 \sqrt{q} \cos \theta $ is an eigenvalue of $G$.

Let 
\[
\phi ( \lambda ) := \left\{
\begin{array}{ll}
\frac{q +1}{2 \pi } \frac{ \sqrt{4q- \lambda {}^2 }}
{(q+1 )^2 - \lambda {}^2 } & \mbox{if $ \mid \lambda \mid \leq 2 \sqrt{q}$, } \\
0 & \mbox{otherwise . }
\end{array}
\right.
\]
Furthermore, the {\em girth} $g(G)$ of a graph $G$ is the minimum of the length 
of prime, reduced cycles in $G$.

\begin{theorem}[Sunada]
Let $ \{ G_n \} {}^{ \infty }_{n=1} $ be a family of 
$(q+1)$-regular graphs such that 
$ \lim_{n \rightarrow \infty } g( G_n )= \infty $. 
For $a,b \in \mathbb{R} (a<b)$, let 
\[
\phi {}_n ([a,b])=| \{ \lambda \in Spec ( G_n ) \mid 
a \leq \lambda \leq b \} | . 
\]
Then 
\[
\lim_{n \rightarrow \infty } \frac{1}{|V( G_n )|}
\phi {}_n ([a,b])= \int^b_a \phi ( \lambda ) d \lambda . 
\]
\end{theorem}

\section{The Grover matrix of a graph and the trace formula with respect to it} 

Let $G$ be a connected graph with $n$ vertices band $m$ edges. 
Then the {\em Grover matrix} ${\bf U} ={\bf U} (G)=( U_{ef} )_{e,f \in D(G)} $ 
of $G$ is defined by 
\[
U_{ef} =\left\{
\begin{array}{ll}
2/d_{t(f)} (=2/d_{o(e)} ) & \mbox{if $t(f)=o(e)$ and $f \neq e^{-1} $, } \\
2/d_{t(f)} -1 & \mbox{if $f= e^{-1} $, } \\
0 & \mbox{otherwise}
\end{array}
\right. 
\]
A quantum walk on $G$ with the Grover matrix as the transition matrix is called a 
{\em Grover walk} on $G$.  

We introduce the {\em positive support} ${\bf F}^+ =( F^+_{ij} )$ of 
a real matrix ${\bf F} =( F_{ij} )$ as follows: 
\[
F^+_{ij} =\left\{
\begin{array}{ll}
1 & \mbox{if $F_{ij} >0$, } \\
0 & \mbox{otherwise}
\end{array}
\right.
\]

If the degree of each vertex of $G$ is not less than 2, i.e., $\delta (G) \geq 2$, 
then $G$ is called a {\em md2 graph}. 
We state a relationship between a discrete-time quantum walk and the Ihara zeta function 
of a graph was given by Ren et al. [15].

\begin{theorem}[Ren, Aleksic, Emms, Wilson and Hancock] 
Let ${\bf B} - {\bf J}_0 $ be the Perron-Frobenius operator (or the edge matrix) 
of a simple md2 graph $G$. 
Let ${\bf U}$ be the Grover matrix of $G$. 
Then the ${\bf B} - {\bf J}_0 $ is the positive support of the transpose of ${\bf U} $, i.e., 
\[
{\bf B} - {\bf J}_0 =( {\bf U}^t )^+ . 
\]
\end{theorem}

Note that 
\[
{\bf Z} (G,u)^{-1} = \det ( {\bf I}_{2m} -u ( {\bf B} - {\bf J}_0 )), \ m=|E(G)|\ (see [5] ) . 
\]

Next, we recall a zeta function on the positive support of 
the Grover matrix of a graph. 
Let $G$ be a connected graph with $n$ vertices and $m$ edges and 
${\bf U} = {\bf U} (G)$ the Grover matrix of $G$. 
By Theorem 4, we obtain the following result.

\newtheorem{proposition}{Proposition} 
\begin{proposition} 
Let $G$ be a connected graph with $m$ edges. 
Then  
\[
{\bf Z} (G, u)^{-1} = \det ( {\bf I}_{2m} -u {\bf U}^+ ) .  
\] 
\end{proposition}

The Ihara zeta function of a graph is just a zeta function on the positive support of 
the Grover matrix of a graph. 

Konno and Sato [9] presented a formula for the characteristic polynomial of ${\bf U} $. 

Let $G$ be a connected graph with $n$ vertices and $m$ edges. 
Then the $n \times n$ matrix ${\bf T} (G)=( T_{uv} )_{u,v \in V(G)}$ is given as follows: 
\[
T_{uv} =\left\{
\begin{array}{ll}
1/( \deg {}_G u)  & \mbox{if $(u,v) \in D(G)$, } \\
0 & \mbox{otherwise. }
\end{array}
\right.
\] 
This is the transition matrix of the simple random walk on $G$.

\begin{theorem}[Konno and Sato]
Let $G$ be a connected graph with $n$ vertices $v_1 , \ldots , v_n$ and $m$ edges. 
Then, for the Grover matrix of $G$,  
\[
\det ( \lambda {\bf I}_{2m} - {\bf U} )
=( \lambda {}^2 -1)^{m-n} \det (( \lambda {}^2 +1) {\bf I}_n -2 \lambda {\bf T} (G)) 
=\frac{( \lambda {}^2 -1)^{m-n} \det (( \lambda {}^2 +1) {\bf D} -2 \lambda {\bf A} (G))}
{d_{v_1} \cdots d_{v_n }} . 
\]
\end{theorem}

Thus, we can express the spectra of the Grover matrix ${\bf U}$ by means of those of 
${\bf T} (G)$(see [3]).

\newtheorem{corollary}{Corollary}
\begin{corollary}[Emms, Hancock, Severini and Wilson] 
Let $G$ be a connected graph with $n$ vertices and $m$ edges. 
The Grover matrix ${\bf U}$ has $2n$ eigenvalues of the form 
\[
\lambda = \lambda {}_T \pm i \sqrt{1- \lambda {}^2_T } , 
\]
where $\lambda {}_T $ is an eigenvalue of the matrix ${\bf T} (G)$. 
The remaining $2(m-n)$ eigenvalues of ${\bf U}$ are $\pm 1$ with equal multiplicities. 
\end{corollary}

Now, we propose a new zeta function of a graph.  
Let $G$ be a connected graph with $m$ edges. 
Then we define a zeta function $ \overline{{\bf Z}} (G, u)$ of $G$ satisfying   
\[
\overline{{\bf Z}} (u)^{-1} =\overline{{\bf Z}} (G, u)^{-1} = \det ( {\bf I}_{2m} -u {\bf U} ) .    
\]

By Theorem 5, we obtain the determinant expression for $\overline{{\bf Z}} (u)$.

\begin{corollary} 
Let $G$ be a connected graph with $n$ vertices $v_1 , \ldots , v_n$ and $m$ edges. 
Then  
\[
\overline{{\bf Z}} (u)^{-1} =(1- u^2 )^{m-n} \det ((1+ u^2) {\bf I}_n -2u {\bf T} (G)) 
=\frac{(1+u^2 )^{m-n} \det ((1+ u^2 ) {\bf D} -2u {\bf A} (G))}{d_{v_1} \cdots d_{v_n }} . 
\]
\end{corollary}

Next, we consider a regular graph. 
Let $G$ be a connected $(q+1)$-regular graph. 
Then we have 
\[
{\bf D} =(q+1) {\bf I}_n . 
\]

The following result is a consequence of Corollary 2.

\begin{theorem}[Konno and Sato]
Let $G$ be a connected $(q+1)$-regular graph with $n$ vertices and $m$ edges. 
Furthermore, set $Spec (G)= \{ \lambda {}_1 , \ldots , \lambda {}_n \} $. 
Then, 
\[
\overline{{\bf Z}} (u)^{-1} =(q+1 )^{-n} (1-u^2 )^{m-n} \det ((q+1)(1+u^2 ) {\bf I}_n -2u {\bf A} (G)) 
\]
\[
=(1-u^2 )^{m-n} \prod^n_{j=1} \left( u^2 - \frac{2 \lambda {}_j }{q+1} u+1 \right) . 
\]
\end{theorem}

Konno, Mitsuhashi, Morita and Sato [10] presented the trace formula for a regular graph $G$.

Now, let $G$ be a connected $(q+1)$-regular graph. 
Furthermore, let $h( \theta )$ be a complex function on $\mathbb{R} $ with period $2 \pi $. 
The Fourier transform of $h( \theta )$ is given by 
\[
\hat{h} (k)= \frac{1}{ 2 \pi } \int^{ 2 \pi }_0 h( \theta ) 
e {}^{ \sqrt{-1} k \theta } d \theta , 
\]
where $k \in \mathbb{Z} $. 
Furthermore, let $\tilde{w} : D(G) \times D(G) \longrightarrow \mathbb{C} $ be defined such that 
\[
\tilde{w} (e,f) =\left\{
\begin{array}{ll}
2/( \deg {}_G t(e))  & \mbox{if $t(e)=o(f)$ and $f \neq e^{-1} $, } \\
2/( \deg {}_G t(e)) -1 & \mbox{if $f= e^{-1} $, } \\
0 & \mbox{otherwise, }
\end{array}
\right.
\] 
and let 
\[
\tilde{w} (C)=\tilde{w} (e_1 , e_2 ) \tilde{w} (e_2 , e_3 ) \cdots \tilde{w} (e_r , e_1 ) , \ C=( e_1 , e_2 , \ldots , e_r ) . 
\]

\begin{theorem}[Konno, Mitsuhashi, Morita and Sato]
Let $G$ be a connected $(q+1)$-regular graph with $n$ vertices and $m$ edges. 
Suppose that $q>1$.
Set $Spec (G)= \{ \lambda {}_1 , \ldots , \lambda {}_n \} $. 
Furthermore, for each $ \lambda {}_j \ (1 \leq j \leq n)$, 
let $z_j = e {}^{ \sqrt{-1} \theta {}_j}$ be a root of the following quadratic equation: 
\[
v^2 - \frac{2 \lambda {}_j }{q+1} v+1=0 ,  
\] 
where $\theta {}_j =0$ if $\lambda {}_j =q+1$. 
Then the following trace formula holds:
\[
\sum^n_{j=1} h( \theta {}_j )= \frac{m}{2} \int^{ \pi }_{0} h( \theta ) \frac{d \theta }{ \pi }
+ \frac{1}{2}\sum_{[C]} \sum^{ \infty }_{k=1} |C| \tilde{w} (C)^k \hat{h} (k|C|) ,  
\] 
where $[C]$ runs over all equivalence classes of prime cycles of $G$. 
\end{theorem}

\section{The twisted Grover matrix of a mixed graph} 

Let $G$ be a connected graph with $m$ edges and $A(G) \subset V(G) \times V(G)$. 
Then the triple $G=(V(G), E(G) , A(G))$ of $V(G), E(G) , A(G)$ is called the {\em mixed graph} (see [11,12]). 
Set $A(G)^{-1} = \{ e^{-1} \mid e \in A(G) \} $. 
We define an $n \times 2m$ matrix ${\bf K} =( K_{ve} )_{v \in V(G), e \in D(G)} $ and 
a $2m \times 2m$ matrix ${\bf C} $ as follows: 
\[
K_{ve} = 
\left\{
\begin{array}{ll}
1/ \sqrt{ \deg v} & \mbox{if $t(e)=v$, } \\
0 & \mbox{otherwise, } 
\end{array}
\right.
\] 
and 
\[
{\bf C} = 2 {\bf K}^* {\bf K} - {\bf I}_{2m} . 
\]
Note that the $(e, f)$-entry $C_{ef} $ of ${\bf C} $ is given by 
\[
C_{ef} = 
\left\{
\begin{array}{ll}
2/ d_{t(e)} & \mbox{if $t(e)=t(f)$ and $f \neq e$, } \\
2/ d_{t(e)} -1 & \mbox{if $f=e$, } \\ 
0 & \mbox{otherwise, } 
\end{array}
\right.
\] 
Furthermore, we have 
\[
{\bf K} {\bf K}^* = {\bf I}_{n} . 
\]

Next, let $\theta : D(G) \longrightarrow \mathbb{C} $ be a function such that 
\[
\theta ( e^{-1} )=- \theta (e) \ for \ each \ e \in A(G)  
\]
and 
\[
\theta (e)=0 \ for \ e \in A(G) \cap A(G)^{-1} . 
\]
Then we define a $2m \times 2m$ matrix ${\bf S}_{\theta } =( (S_{\theta } )_{ef} ) _{e,f \in D(G)} $ as follows: 
\[  
( S_{\theta } )_{ef} = 
\left\{
\begin{array}{ll}
e^{i \theta (f)} & \mbox{if $f= e^{-1} $, } \\
0 & \mbox{otherwise, } 
\end{array}
\right.
\] 
Furthermore, let 
\[
{\bf U}_{\theta } = {\bf S}_{\theta } {\bf C} . 
\] 
In the case of $ \theta \equiv 0$, the matrix ${\bf U}_{\theta } $ is equal to the Grover matrix of $G$. 
The matrix ${\bf U}_{\theta } $ is called the {\em twisted Grover matrix} of $G$. 
Note that 
\[
( {\bf U}_{\theta } )_{ef} =\left\{
\begin{array}{ll}
2 e^{-i \theta (e)} / \deg t(f)  & \mbox{if $t(f)=o(e)$ and $f \neq e^{-1} $, } \\
e^{-i \theta (e)} (2/ \deg t(f) -1) & \mbox{if $f= e^{-1} $, } \\
0 & \mbox{otherwise. }
\end{array}
\right.
\]

Now, let an $n \times n$ matrix ${\bf H}_{\theta } =( H^{( \theta )}_{uv} )_{u,v \in V(G)} $ be given as follows: 
\[  
H^{( \theta )}_{uv} = 
\left\{
\begin{array}{ll}
1 & \mbox{if $(u,v) \in A(G) \cap A(G)^{-1} $, } \\
e^{i \theta (u,v)} & \mbox{if $(u,v) \in A(G) \setminus A(G)^{-1}$, } \\
e^{-i \theta (v,u)} & \mbox{if $(u,v) \in A(G)^{-1} \setminus A(G)$, } \\
0 & \mbox{otherwise, } 
\end{array}
\right.
\] 
Then the matrix ${\bf H}_{\theta } $ is called the {\em generalized Hermitian adjacency matrix} of $G$. 
Furthermore, let 
\[
{\bf D}^{1/2} = 
\left[ 
\begin{array}{ccc}
\sqrt{ \deg v_1 } & \cdots & 0 \\
\ & \ddots & \ \\ 
0 & \cdots & \sqrt{ \deg v_n } 
\end{array} 
\right]
\] 
and 
\[
\tilde{{\bf H}}_{\theta } = {\bf D}^{-1/2} {\bf H}_{\theta } {\bf D}^{-1/2} ,  
\]
where $V(G)= \{ v_1 , \ldots v_n \} $.  

Kubota, Segawa and Taniguchi [11] showed the following result.

\newtheorem{lemma}{Lemma} 
\begin{lemma}
\[
\tilde{{\bf H}}_{\theta } = {\bf K} {\bf S}_{\theta } {\bf K}^* . 
\]
\end{lemma}

Then we obtain a result for the twisted Grover matrix.

\begin{theorem}
Let $G$ be a connected graph with $n$ vertices $v_1 , \ldots v_n $ and $m$ edges, and 
$A(G) \subset V(G) \times V(G)$. 
Furthermore, let $\theta : D(G) \longrightarrow \mathbb{C} $ be a function such that 
\[
\theta ( e^{-1} )=- \theta (e) , \ e \in A(G) \ and \   
\theta (e)=0 , \ e \in A(G) \cap A(G)^{-1} . 
\] 
Then, for the twisted Grover matrix of $G$,  
\[
\det ( {\bf I}_{2m} -u {\bf U}_{\theta } )
=(1- u^2 )^{m-n} \det ((1+ u^2 ) {\bf I}_n -2u \tilde{{\bf H}}_{\theta } ) 
=\frac{(1- u^2 )^{m-n} \det ((1+ u^2 ) {\bf D} -2u {\bf H}_{\theta } )}
{d_{v_1} \cdots d_{v_n }} . 
\]
\end{theorem}

\noindent
{\bf Proof}.  Let $A(G)\cap A(G)^{-1} = \{ e_1 , \ldots , e_p , e^{-1}_1 , \ldots , e^{-1}_p \} $ and 
$A(G) \setminus A(G)^{-1} = \{ f_1 , \ldots , f_q \} $, where $p+q=m$. 
Then we have 
\[
A(G)^{-1} \setminus A(G)=  \{ f^{-1}_1 , \ldots , f^{-1}_q \} . 
\]
Arrange arcs of $G$ as follows: 
\[
e_1 , e^{-1}_1 , \ldots , e_p ,  e^{-1}_p , f_1 ,  f^{-1}_1 , \ldots , f_q , f^{-1}_q . 
\]
Thus, we have 
\[
{\bf S}_{\theta } = 
\left[ 
\begin{array}{ccccccc}
0 & 1 &  &  &  & 0 \\
1 & 0 &  &  &  &   \\ 
 &  & \ddots &  &  &  \\ 
 &  &  & 0 & e^{-i \theta (f_1 )} &  \\
 &  &  & e^{i \theta (f_1 )} & 0 &  \\
0 &  &  &  &  & \ddots 
\end{array} 
\right] 
.
\]

But, if ${\bf A} $ and ${\bf B} $ are an $m \times n$ matrix and an $n \times m$ matrix, then 
\[
\det ( {\bf I}_m - {\bf A} {\bf B} )= \det ( {\bf I}_n - {\bf B} {\bf A} ) . 
\]
Thus, we have 
\[
\begin{array}{rcl}
\  &  & \det ( {\bf I}_{2m} -u {\bf U}_{\theta } )= \det ( {\bf I}_{2m} -u {\bf S}_{\theta } {\bf C} ) \\ 
\  &  &    \\
\  & = & \det ( {\bf I}_{2m} -u {\bf S}_{\theta } (2 {\bf K}^* {\bf K} - {\bf I}_{2m} )) \\ 
\  &  &    \\
\  & = & \det ( {\bf I}_{2m} +u {\bf S}_{\theta } -2u {\bf S}_{\theta } {\bf K}^* {\bf K} ) \\ 
\  &  &    \\
\  & = & \det ( {\bf I}_{2m} -2u {\bf S}_{\theta } {\bf K}^* {\bf K} ( {\bf I}_{2m} +u {\bf S}_{\theta })^{-1} ) 
\det ( {\bf I}_{2m} +u {\bf S}_{\theta } ) \\ 
\  &  &    \\
\  & = & \det ( {\bf I}_{n} -2u {\bf K} ( {\bf I}_{2m} +u {\bf S}_{\theta })^{-1} {\bf S}_{\theta } {\bf K}^* ) 
\det ( {\bf I}_{2m} +u {\bf S}_{\theta } ) . 
\end{array}
\]

Here, 
\[
\begin{array}{rcl}
\  &  & \det ( {\bf I}_{2m} +u {\bf S}_{\theta } ) \\
\  &  &    \\
\  & = & \det ( 
\left[ 
\begin{array}{ccccccc}
1 & u &  &  &  & 0 \\
u & 1 &  &  &  &   \\ 
 &  & \ddots &  &  &  \\ 
 &  &  & 1 & u e^{-i \theta (f_1 )} &  \\
 &  &  & u e^{i \theta (f_1 )} & 1 &  \\
0 &  &  &  &  & \ddots 
\end{array} 
\right] 
) \\ 
\  &  &    \\
\  & = & (1- u^2 )^m .  
\end{array}
\]
Furthermore, we have 
\[
\begin{array}{rcl}
\  &  & ( {\bf I}_{2m} +u {\bf S}_{\theta } )^{-1} \\
\  &  &    \\
\  & = &  
\left[ 
\begin{array}{ccccccc}
1 & u &  &  &  & 0 \\
u & 1 &  &  &  &   \\ 
 &  & \ddots &  &  &  \\ 
 &  &  & 1 & u e^{-i \theta (f_1 )} &  \\
 &  &  & u e^{i \theta (f_1 )} & 1 &  \\
0 &  &  &  &  & \ddots 
\end{array} 
\right]^{-1} \\ 
\  &  &    \\
\  & = & \frac{1}{1- u^2 } 
\left[ 
\begin{array}{ccccccc}
1 & -u &  &  &  & 0 \\
-u & 1 &  &  &  &   \\ 
 &  & \ddots &  &  &  \\ 
 &  &  & 1 & -u e^{-i \theta (f_1 )} &  \\
 &  &  & -u e^{i \theta (f_1 )} & 1 &  \\
0 &  &  &  &  & \ddots 
\end{array} 
\right]  
\\ 
\  &  &    \\
\  & = & \frac{1}{1- u^2 } ( {\bf I}_{2m} -u {\bf S}_{\theta } ) . 
\end{array} 
\]
Thus, 
\[
\begin{array}{rcl}
\  &  & \det ( {\bf I}_{2m} -u {\bf U}_{\theta } ) \\ 
\  &  &    \\
\  & = & (1- u^2 )^m \det ( {\bf I}_{n} -2u/(1- u^2 ) {\bf K} ( {\bf I}_{2m} -u {\bf S}_{\theta }) {\bf S}_{\theta } {\bf K}^* ) \\ 
\  &  &    \\
\  & = &  (1- u^2 )^{m-n} \det ((1- u^2 ) {\bf I}_{n} -2u {\bf K} {\bf S}_{\theta } {\bf K}^* +2 u^2 {\bf K} {\bf S}^2_{\theta } {\bf K}^* ) . 
\end{array}
\]

Since 
\[
{\bf S}^2_{\theta } = {\bf I}_{2m} , 
\]
we have 
\[
{\bf K} {\bf S}^2_{\theta } {\bf K}^* = {\bf K} {\bf K}^* = {\bf I}_{n} . 
\]
By Lemma 1, it follows that 
\[
\begin{array}{rcl}
\  &  & \det ( {\bf I}_{2m} -u {\bf U}_{\theta } ) \\ 
\  &  &    \\
\  & = &  (1- u^2 )^{m-n} \det ((1+ u^2 ) {\bf I}_{n} -2u {\bf K} {\bf S}_{\theta } {\bf K}^* ) \\  
\  &  &    \\
\  & = &  (1- u^2 )^{m-n} \det ((1+ u^2 ) {\bf I}_{n} -2u \tilde{{\bf H}}_{\theta } ) . 
\end{array}
\]

By the definition of $\tilde{{\bf H}}_{\theta } $, we have  
\[
\begin{array}{rcl}
\  &  & \det ( {\bf I}_{2m} -u {\bf U}_{\theta } ) \\ 
\  &  &    \\
\  & = & (1- u^2 )^{m-n} \det ((1+ u^2 ) {\bf I}_{n} -2u {\bf D}^{-1/2} {\bf H}_{\theta } {\bf D}^{-1/2} ) \\  
\  &  &    \\
\  & = & (1- u^2 )^{m-n} \det ( {\bf D}^{-1} ) \det ((1+ u^2 ) {\bf D} -2u {\bf H}_{\theta } ) \\  
\  &  &    \\
\  & = & \frac{(1- u^2 )^{m-n} \det ((1+ u^2 ) {\bf D} -2u {\bf H}_{\theta } )}{d_{v_1} \cdots d_{v_n }} . 
\end{array}
\] 
Q.E.D.

Substituting $u=1/ \lambda $, we obtain the following result.

\begin{corollary}
Let $G$ be a connected graph with $n$ vertices $v_1 , \ldots v_n $ and $m$ edges, and 
$A(G) \subset V(G) \times V(G)$. 
Furthermore, let $\theta : D(G) \longrightarrow \mathbb{C} $ be a function such that 
\[
\theta ( e^{-1} )=- \theta (e) , \ e \in A(G) \ and \   
\theta (e)=0 , \ e \in A(G) \cap A(G)^{-1} . 
\] 
Then, for the twisted Grover matrix of $G$,
\[
\det ( \lambda {\bf I}_{2m} - {\bf U}_{\theta } )
=( \lambda {}^2 -1)^{m-n} \det (( \lambda {}^2 +1) {\bf I}_n -2 \lambda \tilde{{\bf H}}_{\theta } ) 
=\frac{( \lambda {}^2 -1)^{m-n} \det (( \lambda {}^2 +1) {\bf D} -2 \lambda {\bf H}_{\theta } )}
{d_{v_1} \cdots d_{v_n }} . 
\]
\end{corollary}

Thus, we can express the spectra of the twisted Grover matrix ${\bf U}_{\theta } $ by means of those of 
$\tilde{{\bf H}}_{\theta } $.

\begin{corollary}
Let $G$ be a connected graph with $n$ vertices $v_1 , \ldots v_n $ and $m$ edges, and 
$A(G) \subset V(G) \times V(G)$. 
Furthermore, let $\theta : D(G) \longrightarrow {\bf C} $ be a function such that 
\[
\theta ( e^{-1} )=- \theta (e) , \ e \in A(G) \ and \   
\theta (e)=0 , \ e \in A(G) \cap A(G)^{-1} . 
\]
The twisted Grover matrix ${\bf U}_{\theta } $ has $2n$ eigenvalues of the form 
\[
\lambda = \lambda {}_H \pm i \sqrt{1- \lambda {}^2_H } , 
\]
where $\lambda {}_H $ is an eigenvalue of the matrix $\tilde{{\bf H}}_{\theta } $. 
The remaining $2(m-n)$ eigenvalues of ${\bf U}_{\theta } $ are $\pm 1$ with equal multiplicities. 
\end{corollary}

\noindent
{\bf Proof}.  By Corollary 3, we have 
\[
\det ( \lambda {\bf I}_{2m} - {\bf U}_{\theta } )
=( \lambda {}^2 -1)^{m-n} \prod_{ \lambda {}_H \in Spec( \tilde{{\bf H}}_{\theta } )} (( \lambda {}^2 +1) -2 \lambda {}_H \lambda ) . 
\]
Solving $\lambda {}^2 -2 \lambda {}_H \lambda +1=0$, we obtain 
\[
\lambda = \lambda {}_H \pm i \sqrt{1- \lambda {}^2_H } .  
\]
Q.E.D.

Now, we propose a new zeta function of a graph.  
Let $G$ be a connected graph with $m$ edges. 
Then we define a zeta function ${\bf Z}_{\theta } (G, u)$ of $G$ by    
\[
{\bf Z}_{\theta } (u)^{-1} = {\bf Z}_{\theta } (G, u)^{-1} = \det ( {\bf I}_{2m} -u {\bf U}_{\theta } ) .    
\]

By Theorem 8, we obtain the determinant expression for ${\bf Z}_{\theta } (u)$.

\begin{corollary} 
Let $G$ be a connected graph with $n$ vertices $v_1 , \ldots , v_n$ and $m$ edges. 
Then, for the twisted Grover matrix of $G$,  
\[
{\bf Z}_{\theta } (u)^{-1} =(1- u^2 )^{m-n} \det ((1+ u^2) {\bf I}_n -2u \tilde{{\bf H}}_{\theta } ) 
=\frac{(1-u^2 )^{m-n} \det ((1+ u^2 ) {\bf D} -2u {\bf H}_{\theta } )}{d_{v_1} \cdots d_{v_n }} . 
\]
\end{corollary}

Next, we consider a regular graph. 
Let $G$ be a connected $(q+1)$-regular graph. 
Then we have 
\[
{\bf D} =(q+1) {\bf I}_n . 
\]

The following result is obtained from Corollary 5.

\begin{theorem}
Let $G$ be a connected $(q+1)$-regular graph with $n$ vertices and $m$ edges. 
Furthermore, set $Spec ( {\bf H}_{\theta } )= \{ \lambda {}_1 , \ldots , \lambda {}_n \} $. 
Then, 
\[
{\bf Z}_{\theta } (u)^{-1} =(q+1 )^{-n} (1-u^2 )^{m-n} \det ((q+1)(1+u^2 ) {\bf I}_n -2u {\bf H}_{\theta } ) 
\]
\[
=(1-u^2 )^{m-n} \prod^n_{j=1} \left( u^2 - \frac{2 \lambda {}_j }{q+1} u+1 \right) . 
\]
\end{theorem}

By Theorem 9, we obtain the poles for ${\bf Z}_{\theta } (u)$ of a regular graph.

\begin{corollary}
Let $G$ be a connected $(q+1)$-regular graph with $n$ vertices and $m$ edges. 
The zeta function ${\bf Z}_{\theta } (u)$ has $2n$ poles of the form 
\[
\lambda = \frac{ \lambda {}_j \pm \sqrt{ \lambda {}^2_j -(q+1)^2 }}{q+1} , 
\]
where $Spec ( {\bf H}_{\theta } )= \{ \lambda {}_1 , \ldots , \lambda {}_n \} $. 
The remaining $2(m-n)$ poles of ${\bf Z}_{\theta } (u)$ are $\pm 1$ with equal multiplicities. 
\end{corollary}

Now, we estimate the spectrum of the matrix ${\bf H}_{\theta } $.

\begin{proposition} 
Let $G$ be a connected $(q+1)$-regular graph with $n$ vertices and $m$ edges. 
Then, for any eigenvalue $ \lambda $ of ${\bf H}_{\theta } $, we have 
\[
| \lambda | \leq q+1 . 
\]
\end{proposition}

\noindent
{\bf Proof}. Similarly to the proof of Proposition 3.3 in [11], we have 
\[
| \mu | \leq 1 
\]
for each eigenvalue $\mu $ of $\tilde{{\bf H}}_{\theta } $. 
But, we have 
\[
\tilde{{\bf H}}_{\theta } = {\bf D}^{-1/2} {\bf H}_{\theta } {\bf D}^{-1/2} = \frac{1}{q+1} {\bf H}_{\theta } . 
\]
Thus, each eigenvalue $ \lambda $ of ${\bf H}_{\theta } $ is the following form: 
\[
\lambda =(q+1) \mu , \ \mu \in Spec( \tilde{{\bf H}}_{\theta } ) . 
\]
Therefore, it follows that 
\[
| \lambda |=(q+1)| \mu | \leq q+1 . 
\]
Q.E.D.

We use the following result.

\begin{lemma}
Let $G$ be a connected $(q+1)$-regular graph with $n$ vertices and $m$ edges, and 
$Spec ( {\bf H}_{\theta } )= \{ \lambda {}_1 , \ldots , \lambda {}_n \} $.
If $| \lambda {}_j |<q+1\ (1 \leq j \leq n)$, then 
\[
u^2 - \frac{2 \lambda {}_j }{q+1} u+(q+1)=(u- u^+_j )(u- u^-_j ) ,  
\] 
where 
\[ 
u^{ \pm}_j = \frac{ \lambda {}_j \pm \sqrt{ \lambda {}^2_j -(q+1)^2 }}{q+1} . 
\]
\end{lemma}

\noindent
{\bf Proof}. If $| \lambda {}_j |<q+1$, then 
we have 
\[
u^{\pm}_j =x \pm y i, \ x= \frac{ \lambda {}_j }{q+1} , \ y= \frac{ \sqrt{(q+1)^2 - \lambda {}^2_j }}{q+1} , 
\]
i.e., 
\[
x^2 +y^2 =1 
\]
Thus, if $| \lambda {}_j |<q+1$, then 
\[
u^{ \pm}_j = e^{ \pm i \theta {}_j } , i.e., \ u^-_j =1/ u^+_j . 
\]
Therefore, we have 
\[
u^2 - \frac{2 \lambda {}_j }{q+1} u+(q+1)=(u- u^+_j )(u- u^-_j ) . 
\] 
Q.E.D.

Note that, if $ \lambda {}_j =q+1$, then $u_j =1 $.

We present the Euler product and the exponential expression for the new zeta function of 
a graph. 

At first, 
\[
( {}^t {\bf U}_{\theta } )_{ef} =\left\{
\begin{array}{ll}
2 e^{- \theta (f)} / \deg t(e)  & \mbox{if $t(e)=o(f)$ and $f \neq e^{-1} $, } \\
e^{-i \theta (f)} (2/ \deg t(e) -1) & \mbox{if $f= e^{-1} $, } \\
0 & \mbox{otherwise. }
\end{array}
\right.
\] 
Then we give three weight functions $w: D(G) \times D(G) \longrightarrow {\bf C} $ as follows: 
\[
w (e,f) =\left\{
\begin{array}{ll}
2 e^{-i \theta (f)} / \deg t(e)  & \mbox{if $t(e)=o(f)$ and $f \neq e^{-1} $, } \\
e^{-i \theta (f)} (2/ \deg t(e) -1) & \mbox{if $f= e^{-1} $, } \\
0 & \mbox{otherwise. }
\end{array}
\right.
\] 
For a cycle $C=( e_1, e_2 , \cdots , e_r )$, let 
\[
w(C)=w(e_1 , e_2 ) \cdots w( e_{r-1} , e_r) w( e_r , e_1 ) . 
\]

\begin{theorem}
Let $G$ be a connected graph with $n$ vertices $v_1 , \ldots v_n $ and $m$ edges, and 
$A(G) \subset V(G) \times V(G)$. 
Furthermore, let $\theta : D(G) \longrightarrow {\bf C} $ be a function such that 
\[
\theta ( e^{-1} )= \theta (e)^{-1} , \ e \in A(G) \ and \   
\theta (e)=0 , \ e \in A(G) \cap A(G)^{-1} . 
\]
Then the new zeta function of $G$ is given by 
\[
{\bf Z}_{\theta } (u)= \prod_{[C]} (1-w(C) u^{|C|} )^{-1} 
= \exp (\sum^{\infty}_{k=1} \frac{N_k}{k} u^k ) . 
\]
where $[C]$ runs over all equivalence classes of prime cycles in $G$, 
and $N_k $ is defined by 
\[
N_k = \sum \{ w(C) \mid C: \ a \ cycle \ of \ length \ k \ in \ G \} . 
\]
\end{theorem}

\noindent
{\bf Proof}. At first, we have 
\[
{\bf Z}_{\theta } (u)^{-1} = \det ( {\bf I}_{2m} -u {\bf U}_{\theta } ) 
= \det ( {\bf I}_{2m} -u {}^t {\bf U}_{\theta } ) .  
\] 

At first, 
\[
\log {\bf Z}_{\theta } (u)=\log \det ( {\bf I}_{2m} -u {}^t {\bf U}_{\theta } )^{-1} 
=- {\rm Tr} \log ( {\bf I}_{2m} -u {}^t {\bf U}_{\theta } )
= \sum^{\infty}_{k=1} \frac{1}{k} {\rm Tr} [( {}^t {\bf U}_{\theta } )^k ] u^k . 
\]
Since 
\[
w(e,f)=({}^t {\bf U}_{\theta } )_{ef} , \ e,f \in D(G),  
\]
we have 
\[
{\rm Tr} [( {}^t {\bf U}_{\theta } )^k )]= \sum \{ w(C) \mid C: \ a \ cycle \ of \ length \ k \ in \ G \} = N_k . 
\]
Hence, 
\[
\log {\bf Z}_{\theta } (u)= \sum^{\infty}_{k=1} \frac{N_k }{k} u^k . 
\]
Thus, 
\[
{\bf Z}_{\theta } (u)= \exp \left( \sum^{\infty}_{k=1} \frac{N_k }{k} u^k \right) . 
\]

Next, we consider the Euler product of the zeta function of $G$. 
For a cycle $\tilde{C}$ with length $k$, the exists a prime cycle $C$ of length $p$ such that 
$\tilde{C} = C^l $ and $k=pl$. 
Furthermore, note that 
\[
|[\tilde{C}]|=p . 
\]
Thus, we have 
\[
\begin{array}{rcl}
\log \prod_{[C]} (1-w(C) u^{|C|} )^{-1} & = & - \sum_{[C]} \log (1-w(C) u^{|C|} ) \\
\  &  &    \\
\  & = & \sum_{[C]} \sum^{\infty}_{k=1} \frac{w(C)^k }{k} u^{k|C|} , 
\end{array}
\]
where $[C]$ runs over all equivalence classes of prime cycle in $G$. 
Therefore, 
\[
\begin{array}{rcl}
\frac{d}{du} \log \prod_{[C]} (1-w(C) u^{|C|} )^{-1} & = & u^{-1} \sum_{[C]} \sum^{\infty}_{k=1} |C| w(C)^k u^{k|C|} \\
\  &  &    \\
\  & = &  u^{-1} \sum^{\infty}_{k=1} \sum_{[C]} |C| w(C)^k u^{k|C|} \\
\  &  &    \\
\  & = &  u^{-1} \sum^{\infty}_{k=1} \sum_{C} w(C^k ) u^{k|C|} \\
\  &  &    \\
\  & = &  u^{-1} \sum^{\infty}_{k=1} N_k u^k .
\end{array}
\]
Hence, it follows that 
\[
\log \prod_{[C]} (1-w(C) u^{|C|} )^{-1} = \sum^{\infty}_{k=1} \frac{N_k }{k} u^k , 
\]
and so 
\[
\prod_{[C]} (1-w(C) u^{|C|} )^{-1} = \exp ( \sum^{\infty}_{k=1} \frac{N_k }{k} u^k ) , 
\]
Thus, 
\[
{\bf Z}_{\theta } (u)= \prod_{[C]} (1-w(C) u^{|C|} )^{-1} . 
\]
Q.E.D.

\section{The trace formulas for regular graphs}

We present the trace formula for a regular graph $G$.

Now, let $G$ be a connected $(q+1)$-regular graph. 
Furthermore, let $h( \theta )$ be a complex function on $\mathbb{R} $ with period $2 \pi $. 
The Fourier transform of $h( \theta )$ is given by 
\[
\hat{h} (k)= \frac{1}{ 2 \pi } \int^{ 2 \pi }_0 h( \theta ) 
e {}^{ \sqrt{-1} k \theta } d \theta , 
\]
where $k \in \mathbb{Z} $.

\begin{theorem}
Let $G$ be a connected $(q+1)$-regular graph with $n$ vertices and $m$ edges. 
Suppose that $q>1$.
Set $Spec ( {\bf H}_{\theta } )= \{ \lambda {}_1 , \ldots , \lambda {}_n \} $. 
Furthermore, for each $ \lambda {}_j (1 \leq j \leq n)$, 
let $z_j = e {}^{ \sqrt{-1} \theta {}_j}$ be a root of the following quadratic equation: 
\[
v^2 - \frac{2 \lambda {}_j }{q+1} v+1=0 ,  
\] 
where $\theta {}_j =0$ if $\lambda {}_j =q+1$. 
Then the following trace formula holds:
\[
\sum^n_{j=1} h( \theta {}_j )= \frac{m}{2} \int^{ \pi }_{0} h( \theta ) \frac{d \theta }{ \pi }
+ \frac{1}{2}\sum_{[C]} \sum^{ \infty }_{k=1} |C| w(C)^k \hat{h} (k|C|) ,  
\] 
where $[C]$ runs over all equivalence classes of prime cycles of $G$. 
\end{theorem}

\noindent
{\bf Proof}.  The argument is an analogue of 
Venkov and Nikitin's method [19]. 
At first, let 
\[
u_j =u^+_j \ (1 \leq j \leq n) .
\]
Then we have $u^-_j =1/ u_j $. 
By Theorem 9 and Lemma 2, we have 
\[
{\bf Z}_{\theta } (u)^{-1} =(1- u^2 ) {}^{m-n} \prod^n_{j=1} (u^2 - \frac{2 \lambda {}_j }{q+1} u +1) 
=(1- u^2 ) {}^{m-n} \prod^n_{j=1} (u- u_j )(u-1/ u_j ) . 
\] 
Thus, 
\[
- \log {\bf Z}_{\theta } (u)=(m-n) \log (1-u^2 )+ \sum^n_{j=1} \log u(u+1/u- \frac{2 \lambda {}_j }{q+1} ) . 
\]
Therefore, 
\begin{equation}
- \frac{d}{du} \log {\bf Z}_{\theta } (u)= \frac{2(m-n)u}{u^2 -1} + \frac{n}{u} 
+ \sum^n_{j=1} \frac{d}{du} \log \frac{1}{u} (u-u_j )(u- 1/u_j ) . 
\end{equation}

Since $G$ is $(q+1)$-regular, we get  
\[
2m=n(q+1) . 
\]
Thus, we have 
\begin{equation}
\frac{2(m-n)u}{u^2 -1} + \frac{n}{u} = \frac{(2m-n) u^2 -n}{u(u-1)(u+1)} = \frac{n(q u^2 -1)}{u(u-1)(u+1)} . 
\end{equation}

Next, by Theorem 10, we obtain 
\[
{\bf Z}_{\theta } (u)= \prod_{[C]} (1- w(C) u^{|C|} )^{-1} .  
\]
Then we have 
\[
\log {\bf Z}_{\theta } (u)=- \sum_{[C]} \log (1- w(C) u^{|C|} )
= \sum^{ \infty }_{k=1} \frac{1}{k} w(C)^k u^{k|C|} . 
\]
Thus, 
\begin{equation}
\frac{d}{du} \log {\bf Z}_{\theta } (u)= \sum^{ \infty }_{k=1} |C| w(C)^k u^{k|C|-1} .
\end{equation} 
By (1), (2) and (3), it follows that  
\begin{equation}
- \sum^n_{j=1} \frac{d}{du} \log \frac{1}{u}(u- u_j )(u-1/ u_j )= \frac{n(q u^2 -1)}{u(u-1)(u+1)} 
+ \sum_{[C]} \sum^{ \infty }_{k=1} |C| w(C)^k u^{k|C|-1} . 
\end{equation}

Now, if $ u_j $ is not contained in $\mathbb{R} $, then $ \mid u_j \mid =1$, i.e., 
$u_j =e {}^{\sqrt{-1} \theta {}_j } $. 
Moreover, let $h( \theta )$ be an even complex function on $\mathbb{R}$ 
with period $2 \pi $ which is analytically continuable to an 
analytic function over $Im \  \theta < \epsilon ( \epsilon >0)$. 
Furthermore, let $C_1$ be the circle of radius $r<1$ traced in the 
positive direction(counterclockwise), $C_2 $ the circle of radius $R>1$ traced in the 
positive direction, and $K= C_1- C_2 $. 
Furthermore, we consider the following three contour integrals:
\[
Q(h,j)=- \frac{1}{2 \pi i} \oint_{K} h(- \sqrt{-1} \log u) 
\frac{d}{du} \log \frac{1}{u} (u- u_j )(u-1/u_j ) du , 
\]
\[
I(h)= \frac{1}{2 \pi i} \oint_{K} h(- \sqrt{-1} \log u) 
\frac{n(q u^2 -1)}{u(u-1)(u+1)} du , 
\]
\[
H(h,k)= \frac{1}{2 \pi i} \oint_{K} h(- \sqrt{-1} \log u) 
u^{k \mid C \mid -1} du .
\]
Then, by (4), we have 
\begin{equation} 
\sum^n_{j=1} Q(h,j)= I(h)+ \sum^{ \infty }_{k=1} \mid C \mid 
w(C)^k H(h,k) .
\end{equation} 

By the property of the residue theorem, we have 
\[ 
Q(h,j)= \frac{1}{2 \pi i} \oint_{-K} h(- \sqrt{-1} \log u) 
\frac{d}{du} \log \frac{1}{u} (u- u_j )(u-1/u_j ) du=h( \theta {}_j )+h(- \theta {}_j ) 
=2h( \theta {}_j ) 
\]  
since $u_j = e^{\sqrt{-1} \theta {}_j }$ is contained in $-K$. 
Thus, 
\begin{equation} 
\lim_{r \rightarrow 1} Q(h,j)=2h( \theta {}_j ) .
\end{equation}

Furthermore, we have 
\[
H(h,k)= \frac{1}{2 \pi i} \oint_{C_1} h(- \sqrt{-1} \log u) 
u^{k \mid C \mid -1} du . 
\]
Set $u=r e {}^{ \sqrt{-1} \theta } ( 0 \leq \theta \leq 2 \pi )$, 
where $- \log r < \epsilon $.  
Then 
\[  
\begin{array}{rcl} 
H(h,k) & = & \frac{1}{2 \pi i} \int^{2 \pi }_0 h(- \sqrt{-1} \log r e {}^{ \sqrt{-1} \theta }) 
r^{k|C|} e {}^{ \sqrt{-1} k |C| \theta } \frac{\sqrt{-1} r e {}^{ \sqrt{-1} \theta }}{r e {}^{ \sqrt{-1} \theta }} d \theta \\
\  &   &                \\ 
\  & = & 
\frac{1}{2 \pi } \int^{2 \pi }_0 h( \theta -i \log r) r^{k|C|} e {}^{ \sqrt{-1} k |C| \theta } d \theta \\
\  &   &                \\ 
\  & = &  r^{k|C|} \hat{h} (k |C|- i \log r ) . 
\end{array} 
\] 
Thus, we get  
\begin{equation}
\lim_{r \rightarrow 1} H(h,k) = \hat{h} (k |C| ) . 
\end{equation}

Furthermore, we have  
\[
I(h)= \frac{1}{2 \pi i} \oint_{C_1} h(- \sqrt{-1} \log u) 
\frac{n(q u^2 -1)}{u(u-1)(u+1)} du. 
\]
Set $u=r e {}^{ \sqrt{-1} \theta } \ (- \pi \leq \theta \leq  \pi ; r<1)$. 
Then 
\[
I(h)= \frac{n}{2 \pi } \int^{ \pi }_{- \pi } h( \theta -i \log r) 
\frac{q r^2 e^{2i \theta } -1}{(r e^{i \theta } -1)(r e^{i \theta } +1)} d \theta . 
\]
But, 
\[
\lim_{r \rightarrow 1} I(h)= 
\frac{n}{2 \pi } \int^{ \pi }_{- \pi } h( \theta ) \frac{q e^{2i \theta } -1}{( e^{i \theta } -1)( e^{i \theta } +1)} d \theta . 
\]
Thus, we get 
\[
\frac{q e^{2i \theta } -1}{( e^{i \theta } -1)( e^{i \theta } +1)} 
= \frac{(q+1)(1- \cos 2 \theta )-i(q-1) \sin 2 \theta }{2(1- \cos 2 \theta )} 
= \frac{q+1}{2} -i \frac{(q-1) \sin 2 \theta }{2(1- \cos 2 \theta )} . 
\]
Therefore,  
\begin{equation}
\lim_{r \rightarrow 1} I(h)= \frac{n}{ \pi } \{ \int^{ \pi }_{0} h( \theta ) 
\frac{q+1}{2} d \theta - \frac{i}{4} \int^{ \pi }_{- \pi } h( \theta ) 
\frac{(q-1) \sin 2 \theta }{2(1- \cos 2 \theta )} d \theta \}  
= \frac{n(q+1)}{2 \pi } \int^{ \pi }_{0} h( \theta ) d \theta . 
\end{equation} 
By (5), (6), (7) and (8), it follows that 
\[
2 \sum^n_{j=1} h( \theta {}_j )= \frac{n(q+1)}{2 \pi } \int^{ \pi }_{0} h( \theta ) d \theta 
+ \sum_{[C]} \sum^{ \infty }_{k=1} |C| w(C)^k \hat{h} (k|C|) . 
\] 
Hence,  
\[
\sum^n_{j=1} h( \theta {}_j )= \frac{m}{2} \int^{ \pi }_{0} h( \theta )  \frac{d \theta }{ \pi }  
+ \frac{1}{2} \sum_{[C]} \sum^{ \infty }_{k=1} |C| w(C)^k \hat{h} (k|C|) . 
\]
Q.E.D.

\end{document}